\documentclass[12pt]{article}
\usepackage{amsfonts,amsmath}
\usepackage{epsfig,graphics}

\renewcommand{\appendix}{%
\renewcommand{\section}{%
\newpage
\thispagestyle{plain}%
\secdef\Appendix\sAppendix}%
\setcounter{section}{0}%
\renewcommand{\thesection}{\Alph{section}}%
}
\newcommand{\Appendix}[2][?]{%
\refstepcounter{section}%
\addcontentsline{toc}{Addendum}%
{\protect\numberline{\appendixname~\thesection}#1}%
{\flushleft\large\bfseries\appendixname\ \thesection\par#2\par}
\sectionmark{#1}\vspace{\baselineskip}}

\newcommand{\sAppendix}[1]{%
{\flushright\large\bfseries\appendixname\par
\centering#1\par}%
\vspace{\baselineskip}}

\def\be{\begin{equation}}

\def\bea{\begin{eqnarray}}

\def\eea{\end{eqnarray}}

\renewcommand{\theequation}{\thesection.\arabic{equation}}

\normalbaselines \oddsidemargin 0.5cm \evensidemargin 0.5cm
\usepackage{color}

\begin{document}
\pagestyle{empty}

\rightline{}

 \vskip 1cm
\begin{center}

{\huge {\textbf{Confluence of  singularities of differential
equation : A Lie algebras  contraction approach}}}

\vspace{6mm}

{\bf \large   Mohammed Brahim Zahaf $^a$\footnote{E-mail address~:
m\_b\_zahaf@yahoo.fr} and Dominique Manchon $^{b}$\footnote{E-mail
address~: Dominique.Manchon@math.univ-bpclermont.fr}}
 \vspace{2mm}

{\small {\it $^{a}$\footnote{Laboratoire de recherche agr\'e\'e par
le MESRS dans le cadre du Fonds National de la Recherche et du
D\'eveloppement Technologique.}Laboratoire de Physique Quantique de
la Mati\`ere et Mod\'elisations Math\'ematiques (LPQ3M),
\\Centre Universitaire de Mascara, 29000-Mascara, Alg\'erie}}

\smallskip
{\small {\it $^{b}$Laboratoire de Math\'ematiques, CNRS-UMR 6620\\Universit\'e Blaise Pascal\\24 avenue des Landais 63177 Aubi\`ere Cedex,
France }}

\end{center}
\vspace{4cm}
\begin{abstract}
{\small \noindent  We investigate here the confluence of
singularities of Mathieu differential equation by means of the
Lie algebra contraction of the motion group $M(2)$ on the Heisenberg
group $H(3)$.  }
\end{abstract}

\pagestyle{plain} \setcounter{page}{1}


\newpage


%
\renewcommand{\thesection}{\arabic{section}}
\renewcommand{\theequation}{\arabic{section}.{\arabic{equation}}}
\setcounter{equation}{0}
\section{Introduction}
In this paper we deal with second order homogenous differential
equations
\begin{equation}\label{eqdiff}
P_0(z)y''(z)+P_1(z)y'(z)+P_2(z)y(z)=0
\end{equation}
where $P_0(z)$, $P_1(z)$ and $P_2(z)$ are polynomials in the
complex $z$ and  should have no common factors. The singularities
of (\ref{eqdiff}) are defined to be the zeros of the polynomial $
P_0(z)$. In the literature  the understanding of the solutions of
(\ref{eqdiff})  in the neighborhood of singularities is very
crucial.  In general, singularities are branching points of at
least one particular solution of (\ref{eqdiff}). We recall that
there are two types of singularities: the regular (Fuchsian)
singularities $z_k$ type for which $ (z-z_k)P_1(z)/P_0(z)$ and $
(z-z_k)^2P_2(z)/P_0(z)$ are analytical functions in the
neighborhood of $z_k$ and the irregular singularities type which
are not regular. A regular singularity  (with the corresponding
Frobenius solution) is elementary if the difference between two
Frobenius exponents is equal to $\frac{1}{2}$; otherwise, it is
nonelementary \cite{{Ince}}. When two singularities collapse in
one point we get, by the so called {\it confluence}, a new
differential equation with possibly different types of
singularities and with less polynomial
parameters.\\

The principal object of this work consists in studying the confluence of
singularities of Mathieu differential equation to the corresponding
singularities of the harmonic oscillator differential equation. The solutions
of the former will converge to the solutions of the latter, in a sense which
will be precised. We use to that purpose the Lie group contraction procedure
performed from the motion group $M(2)$ towards the Heisenberg group
$H(3)$ (or precisely on the semi-direct product of $H(3)$ with $
\mathbb{Z}_2$). M. Andler and D. Manchon,
attempting in \cite{AM} to develop the pseudo-differential calculus for finite
difference operators, have developed the contraction approach at the
group representation level using the Kirillov orbit method. The formalism of
multiresolution analysis (MRA), developed previously by Mallat
\cite{Mal}, gives them the mean to precise the sense of the limit
transformation on group representations, by transforming the
irreducible representation in the Hilbert space related to the first
group $M(2)$ to the irreducible representation of the Hilbert space
$L^2(\mathbb{R})$ of the second group $H(3)$.\\

 To make the link between
our work and the results of \cite{AM}, we should emphasize that the
solutions of Mathieu and harmonic oscillator differential equations are
eigenvectors of two (unbounded) operators acting
respectively on the Hilbert spaces of the two above
irreducible representations. These operators both come from a specific second order
element in the respective enveloping algebras. The confluence of  singularities is then
interpreted as a contraction procedure.\\

We will use below the notion of s-rank \cite {SLS} which
characterize either regular and irregular singular points. To
introduce them we associate to (\ref{eqdiff}) the {\em symbolic
indicial equation}
\begin{equation}\label{indice}
T(z,D)=0
\end{equation}
with  $T(z,D):=P_0(z)D^2+P_1(z)D+P_2(z)$ where  $D$ is the
differentiation operator and $z$ is a formally independent
variable. The two solutions ($m=1,\,2$) of (\ref{indice}), in the
variable $D$, can be represented in the neighbourhood of finite
singularities $z_k$ by the Puiseux series
\begin{equation}
D_m(z_k)=(z-z_k)^{-\mu_{mk}}\sum^{\infty}_{i=0}{h_{mi}(z-z_k)^{i/2}},\quad
h_{m0}(z_k)\neq 0\end{equation} or by
\begin{equation}
D_m(\infty)=z^{\mu_{m\infty}-2}\sum^{\infty}_{i=0}{h_{mi}(\infty)^{-i/2}},\quad
h_{m0}(\infty)\neq 0
\end{equation}
for the singularity at infinity. The s-rank is then defined
respectively for finite and infinite irregular singularities by
\begin{equation}
R_{z_k}= \max_{m=1,2}(\mu_{mk}),\quad
R_{\infty}=\max_{m=1,2}(\mu_{m\infty}).
\end{equation}
Irregular singular points for half-integer s-rank  are called {\em
ramified} and {\em unramified} for integer s-rank. The s-rank of a regular singularity is defined in a different way, and
turns out to be 1/2 for an elementary regular singularity, and 1 for a
non-elementary one \cite{slav1}. The
set $\{R_{z_1},R_{z_2},...,R_{\infty}\}$ of s-ranks of singular
points of equation  (\ref{eqdiff}) constitutes its s-multisymbol for
which its elements number decreases by one in the case of confluence
of two singularities. The two corresponding $R_{z_k}$'s give rise to
a new one greater than their maximum. If the new one is equal to the
sum of the original ones then the confluence is called strong.
Otherwise, it is called weak.

\section{The harmonic oscillator differential equation and the Heisenberg group $H(3)$}
\setcounter{equation}{0}
The Heisenberg group $H(3)$ can be
introduced by a multiplicative operation defined, on the three
dimensional space $\mathbb{R}\times \mathbb{R}\times \mathbb{R}$, by
\begin{equation}
(a,b,t)(a^{\prime},b^{\prime},t^{\prime})=(a+a^{\prime},b+b^{\prime},t+t^{\prime}+a.b^{\prime}).
\end{equation}
or can be realized as group of upper triangular matrices
\begin{equation}\label{realization1}
\mathbf{h}(a, b, t)=\left(
\begin{array}{ccc}
1 & a & t \\
0 & 1 & b \\
0 & 0 & 1
\end{array}
\right). \end{equation}
 One can easily show that
\begin{equation}
\mathbf{h}(a,b,t)\mathbf{h}(a^{\prime},b^{\prime},t^{\prime})=\mathbf{h}(a+a^{\prime},b+b^{\prime},t+t^{\prime}+a.b^{\prime}).
\end{equation}
The corresponding Lie algebra, noted $h(3)$, is the three
dimensional vector space $V$  generated, in the above realization
(\ref{realization1}),  by the following matrices
\begin{equation}
P=\left(
\begin{array}{ccc}
0 & 1 & 0 \\
0 & 0 & 0 \\
0 & 0 & 0
\end{array}
\right) , Q =\left(
\begin{array}{ccc}
0 & 0 & 0 \\
0 & 0 & 1 \\
0 & 0 & 0
\end{array}
\right) \text {and} \,\ E =\left(
\begin{array}{ccc}
0 & 0 & 1 \\
0 & 0 & 0 \\
0 & 0 & 0
\end{array}
\right)
\end{equation}
which verify the commutation relations
\begin{equation}
\begin{array}{l}
\left[ P,Q\right] =E ,\,\ \left[ E,P\right] = \left[ E,Q\right] =0.
\end{array}
\label{rcom1} \end{equation} Further, the Heisenberg group $H(3)$
admits unitary representations that are introduced on the Hilbert
space $L^2(\mathbb{R})$ of complex functions endowed with the scalar
product
\begin{equation}
\left\langle{f_1,f_2}\right\rangle=\int_{-\infty}^{\infty}f_1(x)\overline
{f_2(x)}dx . \end{equation} For a fixed real $h$, these
representations are defined by the operators
\begin{equation}\label{irrep}
(R^h(\mathbf{h}( a,b,t+ ab/2))f)(x)=e^{ih( t+ ab/2)} e^{ihax}f(x+
b)
\end{equation}
and are irreducible only for $ h \neq 0$. Due to the Stone-von
Neumann theorem, for $h \neq 0$ representations given by
(\ref{irrep}) describe all the irreducible unitary representations
of $H(3)$ whose restriction to the centre is nontrivial, up to
unitary equivalence (see for example \cite{TA,Th}). Furthermore,
there are characters of the form $
\chi_{\alpha,\beta}(a,b,t)=e^{i(\alpha a+\beta b)}$ constituting
representations of dimension one  and
being trivial on the center  of  $H(3)$.\\
  Let us emphasize that the element
$\mathbf{h}(0,0,t)$ belongs to the group center $Z$ of $H(3)$ and
corresponds to the operator
\begin{equation} R^h(\mathbf{h}(0,0,t))=e^{ih t}I. \label{a}
\end{equation}
On the other hand the operator $R^h (\mathbf{h}(a,0,0))$ represents the
multiplication by a function:
\begin{equation} (R^h(\mathbf{h}(a,0,0))f)(x)=e^{ihax}f(x) \label{b}
\end{equation} and $R^h(\mathbf{h}(0,b,0))$ is the shift operator:
 \begin{equation}
(R^h(\mathbf{h}(0,b,0))f)(x)=f(x+ b) \label{mult} \end{equation} In
representation $R^h$, it follows from the formulae
(\ref{a})-(\ref{mult}) that the elements $P$, $Q$ and $E$ of $h(3)$
algebra correspond respectively to the operators $P^h$, $Q^h$ and
$E^h$ given by
\begin{equation}\label{drep}
\begin{array}{c}
(P^h f)(x)=\frac{d}{dx}f(x)\\
(Q^hf)(x)=ihxf(x) \\
(E^h f)(x)=ih f(x).
\end{array}
\end{equation}
and satisfy the same commutation
relations as (\ref{rcom1}). It is possible to combine $P^h$ and
$Q^h$ in one operator:
\begin{equation}
H^{h}=\{P^h\}^2+\{Q^h\}^2
\end{equation}
which represents the Hamiltonian operator of harmonic oscillator
algebra. The associated eigenfunctions $e^h_n(x)$ with eigenvalues
$-\mu=-(2n+1)h$ read
\begin{equation}
e^h_n(x)=(2^{n}n!)^{-{1\over{2}}}(\pi
/h)^{-{1\over{4}}}e^{-h{x^2}/{2}}H_n({ {\sqrt{h}x}}),~~h > 0,~~
n=0,1,2,..., \label{bo1}
 \end{equation}
and form an orthonormal basis of $L^2(\mathbb{R})$; $H_n$ stand for
 Hermite polynomials.
In terms of the differential notation (\ref{drep}) the eigenvalues
equation $H^{h}y=-\mu y$ is nothing else than the harmonic
oscillator differential equation
\begin{equation}\label {oscilharm}
\frac{d^2 y}{dx^2}+(\mu -h^2x^2)y=0
\end{equation}
For our purpose, we use the quadratic transformation $t=x^2$ to put
the later equation (\ref{oscilharm}) in the form
\begin{equation}\label {eqA}
t\frac{d^2 y}{dt^2}+\frac{1}{2}\frac{d y}{dt}+\frac{1}{4}(\mu
-h^2t)y=0.
\end{equation}
This equation admits two singularities: one in $0$ which is
elementary regular and one at $\infty$ which is unramified
irregular; hence it can be characterized by the s-multisymbol
$\{\frac12;2\}$.

\section{The Mathieu differential equation and the motion group $M(2)$}
\setcounter{equation}{0}
 In the canonical
form, the Mathieu equation reads
\begin{eqnarray}\label{mathequ1}
\frac{d^2 y}{ds^2}+\left(a-2q\cos2s\right)y=0
\end{eqnarray}
where $s$, $q$  and  the {\em characteristic value} $a$ belong to $\mathbb{R}$. Among its solutions there are the {\em
pseudo-periodic\/} ones (called Floquet solutions ) of the form
\cite{Meixner, AS}
\begin{eqnarray}
me_\lambda(s,~q):=e^{\lambda i
s}\sum_{k=-\infty}^{+\infty}C_{k}^{\lambda}e^{2kis}
\end{eqnarray}
where the coefficients $C_{k}^{\lambda}$ satisfy the recursion
relation
\begin{eqnarray}
(a-(2k+\lambda)^2)C_{k}^{\lambda}-q(C_{k+2}^{\lambda}+C_{k-2}^{\lambda})=0.
\end{eqnarray}
Obviously we have
\begin{eqnarray}
me_\lambda(s+\pi,~q)=e^{\lambda\pi i}me_\lambda(s,~q)
\end{eqnarray}
and
\begin{eqnarray}
me_{-\lambda}(s,~q)=me_\lambda(-s,~q).
\end{eqnarray}
Of particular interest in physics and mathematics is the case
$e^{\lambda\pi i}=\pm1$, so that there exists at least one periodic
solution of period $\pi$  or $2\pi$. In this case and when
$q\rightarrow +\infty$ the characteristic value $a$ can be
approximated by
\begin{eqnarray}\label{valcarac}
a=-2q+2(2n+1)q^{1/2}-\frac{(2n+1)^2+1}{8}+O(n^3q^{-1/2}),~n\in
\mathbb{N}
\end{eqnarray}
and for each $q$ the periodic solution is either even or odd
(often denoted by $ce_n(s, q)$ or $se_{n+1}(s, q)$). The change of
variable $x=\cos^2s$ in the Mathieu equation (\ref{mathequ1})
leads to its algebraic form
\begin{equation}\label{Mathieualg1}
x(1-x)\frac{d^2 y}{dx^2}+\frac{1}{2}(1-2 x)
\frac{dy}{dx}+\frac{1}{4}(a+2q -4qx)y=0
\end{equation}
which admits two elementary regular singular points at $0$ and $1$
and a ramified irregular singularity at infinity, so its
s-multisymbol
is $\{\frac12;\frac12;\frac{3}{2}\}$.\\

We will now consider the Lie algebra, denoted by $m_\alpha(2)$,
generated by the three  generators  $P_\alpha$, $Q_\alpha$ and
$E_\alpha$ with the commutation relations \cite{AM}
\begin{equation}
\begin{array}{l}
\left[ P_\alpha,Q_\alpha\right] =E_\alpha ,\,\ \left[
P_\alpha,E_\alpha\right] =-\alpha^2Q_\alpha ,\,\ \left[
E_\alpha,Q_\alpha\right] =0
\end{array}
\label{rcom4} \end{equation} For $\alpha\neq0$, this Lie algebra
corresponds to the group $\tilde G_\alpha=\mathbb{R}\times
\mathbb{R}^2$ equipped with the semi-direct product
\begin{equation}
\left(\theta, v\right).\left(\theta^\prime,
v^\prime\right)=\left(\theta+\theta^\prime,
v+k_\alpha\left(\theta\right).v^\prime\right) \end{equation} where
\begin{eqnarray} k_\alpha\left(\theta\right)=\left(
\begin{array}{cc}
\cos\left(\alpha\theta\right) & -\alpha\sin\left(\alpha\theta\right)\\
\alpha^{-1}\sin\left(\alpha\theta\right) &
\cos\left(\alpha\theta\right).
\end{array}
\right) \end{eqnarray} The set of all $k_\alpha(\theta)$, for
$\theta \in \mathbb{R}$, is denoted by  $SO_\alpha(2)$. The group
$\tilde G_\alpha$  is the simply connected covering of the group
$M_\alpha(2)=SO_\alpha(2)\times\mathbb{R}^2\equiv
\mathbb{R}/_{2\pi\alpha^{-1}\mathbb{Z}}\times \mathbb{R}^2$ with the
composition law
\begin{eqnarray}
\left(\dot\theta, v\right).\left(\dot\theta^\prime,
v^\prime\right)=\left(\theta\dot+\theta^\prime,
v+k_\alpha\left(\theta\right).v^\prime\right)
\end{eqnarray}
where $\dot \theta$  designs the equivalence class of  $\theta$ and
$k_\alpha(\dot\theta)=k_\alpha\left(\theta\right)$. For $\alpha=1$,
$M_\alpha(2)$  is the euclidian motion group of the plane $M(2)$.
For arbitrary $\alpha$,  $M_\alpha(2)$  is the group of
displacements associated with the euclidean structure defined on
$\mathbb{R}^2$ by
\begin{equation}
{\parallel\left(v_1, v_2 \right)
\parallel}_\alpha^2=v_1^2+\alpha^2v_2^2.
\end{equation}
This is why we call it  the elliptic motion group of plane.  For
$\lambda \in \mathbb{R}\slash\mathbb{Z}$ we introduce the Hilbert
space ${\cal{H}} ^{\alpha,\lambda}$ of functions over $\mathbb{R}$
such that $f(\psi+2\pi k\alpha^{-1})=e^{2i\pi k\lambda}f(\psi)$
and which are square integrable over the pseudo-period
$[0,2\pi\alpha^{-1}]$. On ${\cal{H}} ^{\alpha,\lambda}$ and for
the real  $h\ne 0$, we define an unitary irreducible
representation of $\tilde G_\alpha$ by
\begin{equation}
\left(R^{\alpha, \lambda}_{h}(g)
f\right)(\psi)=e^{ih\left(v_2\cos(\alpha \psi)+\alpha^{-1}
v_1\sin(\alpha\psi)\right)}f(\psi+\theta).
 \end{equation}
where $g(\theta,v)=g(\theta,v_1,v_2) \in \tilde G_\alpha$ and
$f\in {\cal{H}} ^{\alpha,\lambda}$. This representation factorizes
in a representation of $M_\alpha(2)$  if and only if $f$ is
$2\pi\alpha^{-1}$-periodic i.e. $\lambda=0$. The infinitesimal
operators are then:
\begin{eqnarray}
(P^{\alpha, \lambda}_{h}{f})(\psi):=(R^{\alpha, \lambda}_{h}(P){f})(\psi)&=&\frac{\partial }{\partial \theta}(R^{\alpha, \lambda}_{h}(g(\theta,0,0)){f})(\psi) |_{\theta=0}=\frac{d{f}}{d\psi}(\psi)\\
(Q^{\alpha, \lambda}_{h}{f})(\psi):=(R^{\alpha, \lambda}_{h}(Q){f})(\psi)&=&\frac{\partial }{\partial v_1}(R^{\alpha, \lambda}_{h}(g(0,v_1,0)){f})(\psi) |_{v_1=0}=ih\alpha^{-1}\sin(\alpha \psi)f(\psi)\nonumber\\
(E^{\alpha, \lambda}_{h}{f})(\psi):=(R^{\alpha,
\lambda}_{h}(E){f})(\psi)&=&\frac{\partial}{\partial
v_2}(R^{\alpha,
\lambda}_{h}(g(0,0,v_2)){f})(\psi)|_{v_2=0}=ih\cos(\alpha
\psi)f(\psi)\nonumber
\end{eqnarray}
We check easily that the operators $P^{\alpha, \lambda}_{h}$,
$Q^{\alpha, \lambda}_{h}$ and $E^{\alpha, \lambda}_{h}$ verify the
same relations as (\ref{rcom4}), and tend formally to the operators in (\ref{drep}) as $\alpha \rightarrow 0$.\\

On the other hand, we know  that the motion group $M(2)$  is the
symmetry group of the Helmholtz equation
$[\bigtriangleup_2+\omega]\psi(x,y)=0$, where $\bigtriangleup_2
=\frac{\partial^2}{\partial x^2}+\frac{\partial^2}{\partial y^2}$.
In \cite {MI}, Miller had shown that the resolution of the above
equation by the method of separation of variables could be realized
in four orthogonal systems of coordinates: cartesian, polar,
parabolic and elliptic associated respectively to four symmetric
quadratic operators $L_1=Q^2$, $L_2=P^2$, $L_3=\{P,Q\}$ and
$L_4=P^2+d^2Q^2$ in the enveloping algebra of $M(2)$. Here $P$ and
$Q$ are respectively the infinitesimal rotation and translation. $L$
is called symmetric operator of Helmholtz equation if $[L,
\verb"Q"]=R(x,y)\verb"Q"$ with $\verb"Q"= \bigtriangleup_2+\omega$
and  $R(x,y)$ is a complex function defined on $\mathbb{R}^2$.
Corresponding to the representation of $M(2)$ every operator $L_i$
corresponds naturally to a symmetric operator  on a domain of
definition $D $ in ${\cal H}^{\alpha,0}\sim L^2(S^1)$  and each of
them
can be extended to a self-adjoint operator defined on a domain $D'\supseteq D$.\\

In the case of $\tilde G_\alpha $, where the above situation is very
similar, and on the space of $C^{2}$ functions in ${\cal
H}^{\alpha,\lambda}$  the elliptic operator
$L^{\alpha}_4=\left\{P^{\alpha,
\lambda}_{h}\right\}^2+d^2\left\{Q^{\alpha, \lambda}_{h}\right\}^2$,
associated to elliptic coordinate system, corresponds (for $d=1$) to
\begin{equation}\label{operat}
L^{\alpha}_4=\frac{d^2}{d\psi^2}-h^2\alpha^{-2}\sin^2(\alpha
\psi).
\end{equation}
We have dropped the superscript $\lambda$, which does not appear in the
right-hand side of equation \eqref{operat}. The operator $L^{\alpha}_4$
depends on it through the Hilbert space ${\cal H}^{\alpha,\lambda}$ on
which it acts. The equation
\begin{equation}\label{mathieu1}
L^{\alpha}_4 y=-\mu y
\end{equation}
after the change of variable $s=\alpha\psi+\frac{\pi}{2}$, is
nothing else than the Mathieu equation with
\begin{equation}\label{par}
a=\alpha^{-2}\mu-2q ~~\text {and}~~    q=\frac{h^2 \alpha^{-4}}{4}.
\end{equation}
 Setting $t=\alpha^{-2}\sin^2(\alpha\psi)$ in
(\ref{mathieu1}), we find the {\em deformed }algebraic form of the
Mathieu equation
\begin{equation}\label{Mathieualg2}
t(1-\alpha^2 t)\frac{d^2 y}{dt^2}+\frac{1}{2}\lbrace 1-2\alpha^2 t
\rbrace \frac{dy}{dt}+\frac{1}{4}(\mu -h^2t)y=0
\end{equation}
This equation admits three singular points: $0$, $\alpha^{-2}$ which
are elementary regular and $\infty$  which is ramified irregular.
Formally when $\alpha \rightarrow 0$ this equation tends to the
harmonic oscillator differential equation (\ref{eqA}). It is exactly
this strong confluence that will be interpreted in terms of Lie
algebra contraction in the  next section.
\section{The contraction of $M(2)$ on $H(3)$ and the confluence }
\setcounter{equation}{0}
 Let consider the vector space $V$
underlying to the Heisenberg algebra $h(3)$. It is generated by the
basis $P$, $Q$ and  $E$. We denote by $\left[.,.\right]_0 $ the Lie
bracket:
\begin{equation}
\begin{array}{l}
\left[ P,Q\right]_0 =E ,\,\ \left[ P,E\right]_0 = \left[
E,Q\right]_0 =0,
\end{array}
\label{rcom2} \end{equation}
so that $V$ endowed with this bracket
$\left[.,.\right]_0 $ is isomorphic to $h(3)$. We also denote by
$\left[.,.\right]_1 $ the Lie bracket defined by
\begin{equation}
\begin{array}{l}
\left[ P,Q\right]_1 =E ,\,\ \left[ P,E\right]_1 =-Q ,\,\ \left[
E,Q\right]_1 =0,
\end{array}
\label{rcom3} \end{equation} so that the vector space  $V$ equipped
with $\left[.,.\right]_1 $ is isomorphic to the Lie algebra $m(2)$. Let
us introduce the following automorphism $\Phi _\alpha$ of $V$:
\begin{equation}
\Phi _\alpha\left(P\right)=\alpha P,\,\,\,\,\Phi
_\alpha\left(Q\right)=\alpha Q ,\,\,\,\,\Phi
_\alpha\left(E\right)=\alpha^2 E
\end{equation}
and the Lie  bracket $\left[.,.\right]_\alpha $ defined by
\begin{equation}
\left[X,Y\right]_\alpha=\Phi_\alpha^{-1}
\left(\left[\Phi_\alpha\left(X\right),\Phi_\alpha\left(Y\right)\right]_1\right)
. \label{cr1}
\end{equation} Then we have
\begin{equation}
\begin{array}{l}
\left[ P,Q\right]_\alpha =E ,\,\ \left[ P,E\right]_\alpha
=-\alpha^2Q ,\,\ \left[ E,Q\right]_\alpha =0.
\end{array}
\end{equation}
It is obvious that $V$ equipped with the Lie bracket
$\left[.,.\right]_\alpha $ is isomorphic to $m_\alpha(2)$ and
\begin{equation}
\lim_{\alpha\rightarrow 0}\left[X,Y\right]_\alpha
=\left[X,Y\right]_0. \end{equation} This means that the algebra
$h(3)$ is a contraction of $m(2)$.\\
 The authors of \cite{AM} have shown,
using the orbits method, that when $\alpha$ goes to $0$ the group
$\tilde G_\alpha$ "tends" to a degree two extension $G_0$ of
the Heisenberg group; and that the representation $R^{\alpha,
\lambda}_{h}$, acting on ${\cal {H}} ^{\alpha, \lambda}$, converges
topologically (in the sense of Fell) to an irreducible representation of
$G_0$ which restricts on $H(3)$
to the direct sum $R^{h}\oplus R^{-h}$ which acts on
$L^2(\mathbb{R})\oplus L^2(\mathbb{R})$. The geometric picture of this
phenomenon is the following: the
coadjoint orbits of $\tilde G_\alpha$ are cylinders of elliptic
base, converging to the union of two planes of height $ \pm h$ when $\alpha\to
0$, i.e. when the big axis of the ellipse grows to infinity. These two planes
together form a coadjoint orbit of $G_0$.\\

To give a sense to the limit of representations, all the spaces
${\cal{H^{\alpha,\lambda}}}$ should be compared together and also
with $L^2(\mathbb{R})\oplus L^2(\mathbb{R})$, using the
multiresolution analysis according to Mallat and Meyer \cite{Mal,
Me, Db}, see
appendix 1.\\

Finally the confluence holds when $\alpha\rightarrow 0$; the
elliptic operator $L^{\alpha}_4$ tends formally to the Hamiltonian
operator $H^h$ of oscillator harmonic algebra, and in virtue of
(\ref{valcarac}) and (\ref{par}) $\mu$ tends to $(2n+1)h$ and the
equation (\ref{Mathieualg2}) becomes the equation (\ref{eqA}). At
the same time  the differential equation solutions experience the
following limits (recall that $h$ is a fixed positive parameter, and that $q$
and $\alpha$ are related by the equality $q=\frac 14 h^2\alpha^{-4}$):
\begin{eqnarray}
\lim_{\alpha\rightarrow
  0}\frac{ce_{2n}(\alpha\psi+\frac{\pi}{2},q)}{ce_{2n}(\frac{\pi}{2},q)} &=&
\frac{\Gamma(\frac{1}{2}-n)}{2^n\pi^{\frac{1}{2}}}~e^{-h{\psi^2}/{2}}H_{2n}({
  {\sqrt{h}\psi}})\nonumber\\
&=&\frac{(-1)^n2^nn!}{(2n)!}e^{-h{\psi^2}/{2}}H_{2n}({
  {\sqrt{h}\psi}})\nonumber\\
\lim_{\alpha\rightarrow
0}\frac{se_{2n+2}(\alpha\psi+\frac{\pi}{2},q)}{se'_{2n+2}(\frac{\pi}{2},q)}
&=&\frac{\Gamma(-\frac{1}{2}-n)}{2^{n+1}\pi^{\frac{1}{2}}}~
e^{-h{\psi^2}/{2}}H_{2n+1}({ {\sqrt{h}\psi}})\nonumber\\
&=&\frac{(-1)^{n+1}2^{n+1}(n+1)!}{(2n+2)!}e^{-h{\psi^2}/{2}}H_{2n+1}({
  {\sqrt{h}\psi}}).
\end{eqnarray}
 In fact and according to \cite{Meixner} (Satz 10 paragraph 2.333), we have in the interval
$[0,\pi]$
\begin{eqnarray}
\left. \begin{array}{c}
 ce_{n}(z,q)  \\
 se_{n+1}(z,q)
  \end{array}
  \right\}=\left(\frac{\pi q^{\frac12}}{2}\right)^{\frac14}(n!)^{-\frac12}D_{n}(2q^{\frac14}\cos
  z)+O(q^{-\frac38})
\end{eqnarray} as $q\rightarrow\infty$,
 where $D_{m}(\zeta)$ is the parabolic cylinder
function given by
$$D_{m}(\zeta)=\frac{1}{2^{m/2}} e^{-\frac{\zeta^2}{4}}H_m\left(\frac{\zeta}{\sqrt{2}}\right).$$
\section{Conclusion and discussion}
In this work we have examined the singularities confluence of the
Mathieu differential equation towards the harmonic oscillator
differential equation via the Lie algebra contraction of the motion
Lie algebra $m(2)$ to the Heisenberg Lie algebra $h(3)$. The use of
the contraction method to interpret successfully the confluence was
based on the approach of the work \cite{AM}. Someone can now be
tempted to develop  similar interpretations of singularities
confluence for other differential equations. Among many examples we
can cite the case of Lam\'e differential equation with four regular
singularities $0$, $1$, $a$ (which are elementary) and $\infty$
associated with the s-multisymbol $\{\frac12;\frac12;\frac12;1\}$

\begin{equation}\label{lame}
\frac{d^2 y}{dx^2}+\frac{1}{2}\left[
\frac{1}{x}+\frac{1}{x-1}+\frac{1}{x-a} \right]
\frac{dy}{dx}+\frac{\mu-l(l+1)x}{4x(x-1)(x-a)}y=0
\end{equation}
If $a=r^{-2}$ and $x=sn^2(z,r)$, where $sn(z,r)$ is the Jacobi
elliptic function, this equation becomes the Lam\'e differential
equation in the jacobian form
\begin{eqnarray}
\frac{d^2 y}{dz^2}-r^2\left\{l(l+1)sn^2(z,r)-\mu\right\}y=0
\end{eqnarray}
where we have used the fact
\begin{eqnarray}
\frac{d}{dz}sn(z,r)=cn(z,r)dn(z,r)=\sqrt{(1-sn^2(z,r))(1-r^2sn^2(z,r))}
\end{eqnarray}
and so
\begin{eqnarray}
\frac{dx}{dz}=2sn(z,r)cn(z,r)dn(z,r)=2\sqrt{x(1-x)(1-r^2x)}
\end{eqnarray}
 The equation (\ref {lame}), which is a special case of Heun
differential equation, is related to the group $SO_0(2,1)$. In
fact it was shown in \cite{W} that the resolution by separation of
variables of the Laplacian equation $\verb"Q" f=l(l+1)f$ on the
hyperboloid $x_0^2-x_1^2-x_2^2=1$,
 where $\verb"Q"=K_1^2+K_2^2-M_3^2$ and $K_1=-x_0\partial
_{x_2}-x_2\partial _{x_0}$,
$K_2=-x_0\partial_{x_1}-x_1\partial_{x_0}$ and
$M_3=x_1\partial_{x_2}-x_2\partial_{x_1}$, can be realized in nine
coordinate systems associated with nine symmetric quadratic
operators in the enveloping algebra of $SO_0(2,1)$. The above
operators correspond to symmetric operators on the domain $D$ of
$C^{\infty}$ functions in ${\cal H}=L^2(S^1)$ corresponding to the
principal series representations of $SO_0(2,1)$ and each operator
can be extended to one or more self-adjoint operators on ${\cal H}$
(see ref. \cite {KM}). Of special interest is the elliptic operator
$L_{E}=M_3^2+k^2K_2^2$ ($k\in \mathbb{R}$) associated to the
elliptic coordinate system. It corresponds on the Hilbert space
${\cal H}=L^2(S^1)$ to
$$L_E=(1+k^2\cos^2\theta)\frac{d^2
}{d\theta^2}+k^2(2l-1)\sin\theta\cos\theta\frac{d}{d\theta}+k^2(l^2\sin^2\theta+l\cos^2\theta)$$
corresponding to the principal series
($l=-\frac{1}{2}+i\rho,~~0<\rho<\infty$). To retrieve Lam\'e
equation (\ref{lame}) we use the variable change
\begin{eqnarray}
f(\theta)=(1+k^2\cos^2\theta)^{l/2}y(x),
\end{eqnarray}
and $x=\frac{\sin^2\theta}{1+k^2\cos^2\theta}$, $a=-\frac{1}{k^2}$
in the equation $L_E f=\mu k^2 f$.

 On the other hand, the contraction of the group
$SO_0(2,1)$ on the Heisenberg group $H(3)$ brings the possibility
to interpret the double confluence of the Lam\'e differential
equation to harmonic oscillator differential equation. In fact, we
set
\begin{eqnarray}
P^{\alpha}=\alpha M_3,~~Q^{\alpha}=-\alpha^2K_2,~\text{and}
~E=\alpha^3K_1
\end{eqnarray}then we have
\begin{eqnarray}
[P^{\alpha},Q^{\alpha}]=E^{\alpha},~~[E^{\alpha},P^{\alpha}]=\alpha^2Q^{\alpha},~~[E^{\alpha},Q^{\alpha}]=\alpha^4P^{\alpha}
\end{eqnarray}
and in the limit $\alpha\rightarrow 0$ we obtain the commutation
relations (\ref{rcom1}) of the Heisenberg Lie algebra. Similarly as
above, setting
$f(\theta)=(1+\alpha^2k^2\cos^2\alpha\theta)^{l/2}y(x)$ and
$x=\frac{\alpha^{-2}\sin^2\alpha\theta}{1+\alpha^2k^2\cos^2\alpha\theta}$
in the equation $L^{\alpha}_{E}f=-\mu  f$, where
$L^{\alpha}_{E}=\left\{P^{\alpha}\right\}^2+k^2\left\{Q^{\alpha}\right\}^2$
($k\in \mathbb{R}$) is the deformed elliptic operator associated to
the elliptic coordinate system and corresponding to the principal
series on $D$ in the Hilbert space ${\cal H}=L^2(S^1)$, we obtain
the deformed Lam\'e equation
\begin{equation}\label{deformedlame}
\frac{d^2 y}{dx^2}+\frac{1}{2}\left[
\frac{1}{x}+\frac{\alpha^2}{\alpha^2x-1}+\frac{\alpha^4k^2}{\alpha^4k^2x+1}
\right] \frac{dy}{dx}+\frac{-\mu -\alpha^6k
^2l(l+1)x}{4x(\alpha^2x-1)(\alpha^4k^2x+1)}y=0
\end{equation}
In the limit $\alpha\rightarrow 0$, after setting
$\rho=\alpha^{-3}k^{-1}h$, this  equation becomes the harmonic
oscillator differential equation (\ref{eqA}). In this process, the
two regular singularities $1$ and $a$ coalesce to $\infty$ and this
double confluence is also strong. Furthermore the contraction of
$SO_0(2,1)$ on the motion group $M(2)$ \cite{Meck} permits us to
interpret the confluence of the Lam\'e differential equation to the
Mathieu differential equation. In this case the elementary regular
singular point $a$ coalesces to $\infty$. Moreover, the result that
the periodic solutions of the Lam\'e equation tend to periodic
solutions of Mathieu equation already exists in the work of Kalnins
{\em et al.} \cite{KMP} given in terms of the contraction of $SO(3)$
to $M(2)$. Therein,  the irreducible representation of $SO(3)$ is
labeled by $l$ (integer)
 and the periodic Lam\'e solutions are polynomials
which experience the following limits
\begin{eqnarray}
&&\lim_{l\rightarrow
+\infty}\prod_{j=1}^{2m}(1-\frac{u}{\theta_j})=\frac{ce_{2m}(\theta,q)}{ce_{2m}(0,q)},\quad\lim_{l\rightarrow
+\infty}u^{\frac12}\prod_{j=1}^{2m+1}(1-\frac{u}{\theta_j})=\frac{se_{2m+1}(\theta,q)}{se'_{2m+1}(0,q)}\\
&&\lim_{l\rightarrow
+\infty}(1-u)^{\frac12}\prod_{j=1}^{2m+1}(1-\frac{u}{\theta_j})=\frac{ce_{2m+1}(\theta,q)}{ce_{2m+1}(0,q)},\quad
\lim_{l\rightarrow
+\infty}u^{\frac12}(1-u)^{\frac12}\prod_{j=1}^{2m+2}(1-\frac{u}{\theta_j})=\frac{se_{2m+2}(\theta,q)}{se'_{2m+2}(0,q)}\nonumber
\end{eqnarray}
where $u=\sin^2\theta$, $q=\frac{h^2}{4}$, $a=\frac{l^2}{h^2}$ and
the zeros $\theta_j$ satisfy  some fundamental relation which
depends on $a$. Finally we can summarize this study by the following
commutative diagram which exhibits the interpretation of confluence
phenomena by the Lie algebra contraction:

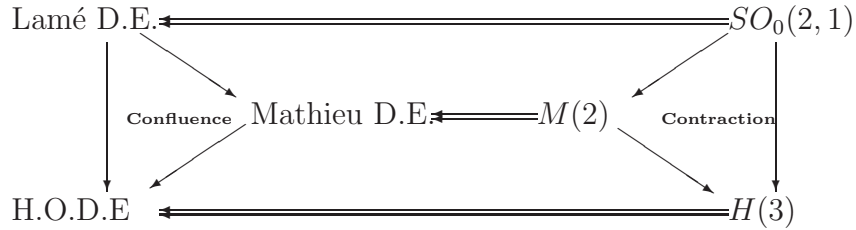
\setlength{\unitlength}{.5in}
\begin{center}
\begin{picture}(8,4)
\put(-0.5,3){Lam\'e D.E.} \put(2,2){Mathieu D.E.}
\put(-0.5,1){H.O.D.E } \put(7,3){$SO_0(2,1)$} \put(5,2){$M(2)$}
\put(7,1){$H(3)$} \put(0.5,2.85){\vector(0,-1){1.55}}
\put(7.5,2.85){\vector(0,-1){1.55}} \put(7,3.1){\vector(-1,0){5.95}}
\put(7,1.1){\vector(-1,0){5.95}} \put(0.85,2.95){\vector(3,-2){1}}
\put(5.85,1.95){\vector(3,-2){1}} \put(7,2.95){\vector(-3,-2){1}}
\put(1.95,2){\vector(-3,-2){1}} \put(5,2.1){\vector(-1,0){1.10}}
\put(7,3.05){\vector(-1,0){5.95}} \put(7,1.05){\vector(-1,0){5.95}}
\put(5,2.05){\vector(-1,0){1.10}} \put(0.7,2){{\bf\tiny Confluence}}
\put(6.3,2){{\tiny \bf Contraction}}
\end{picture}
\end{center}
\vskip 0.5cm \noindent{\bf Acknowledgements:} {\em One of us (M. B.
Z.) would like to thank Prof. H. Dib and Dr. A. Yanallah for
precious help and useful discussions.
 }

\newpage
\setcounter{section}{0} \Appendix[0]{Multiresolution analysis} In
this appendix we recall some background on Multiresolution analysis.
So it is useful in below  to identify ${\cal{H}}^{\alpha,\lambda}$
space with ${\cal {V}}^{\alpha,\lambda}=\ell^{2}(\alpha
(\mathbb{Z}+\lambda))$ space by the isometric isomorphism ${\cal
F}:f\mapsto \hat f $ defined, for $x\in \alpha (\mathbb{Z}+\lambda)$
, by
\begin{equation}
\hat f(x)=\int_{\mathbb{R}/2\pi
\alpha^{-1}\mathbb{Z}}e^{-ix\psi}f(\psi)d\psi.
\end{equation}
We recall briefly the construction of  Littlewood-Paley-Meyer (LPM)
wavelets \cite{AM, Me}. Let $\phi$ a function of class $C^{\infty}$
such that \footnote{As an example of $\phi$ function we can adopt
the one of \cite{YO}:\\$$\phi(\psi)=\sqrt{g(\psi)g(-\psi)}$$ where
$g(\psi)=\frac{h(4\pi/3-\psi)}{h(\psi-2\pi/3)+h(4\pi/3-\psi)}$,
$h(\psi)=\exp(-1/ \psi^2)$, $\psi>0$}
\begin{equation}
\begin{array}{rl}
 \phi(\xi)=1 & \qquad \forall~|\xi|\leq {2\pi\over 3}, \\
0< \phi(\xi)<1 & \qquad \forall~{2\pi\over3}<|\xi|\leq{4\pi\over3},\\
 \phi(\xi) =0&\qquad\forall~|\xi|\geq{4\pi\over3},\\
\phi^2(\xi)+\phi^2(2\pi-\xi)=1&\qquad \forall ~0\leq\xi\leq 2\pi,
\end{array}
\end{equation}
We denote by $V^{\alpha,\lambda}$ the subspace of $L^2(\mathbb{R})$
of functions  having the Fourier transformation of the form
\begin{equation}
m(\xi)\phi(\alpha\xi)
\end{equation}
where $m \in {\cal{H^{\alpha,\lambda}}}$ .\\
Using the proprieties of function $\phi$, we remark that
\begin{equation}\label{per}
\sum_{k\in \mathbb{Z}}\phi^2(\xi+2k\pi)=1
\end{equation}
which leads to
\begin{equation}
\int_{-\infty}^{+\infty}m_1(\xi)\overline{m_2(\xi)}\phi^2(\alpha\xi)d\xi=
\int_{0}^{2\pi\alpha^{-1}}m_1(\xi)\overline{m_2(\xi)} d\xi
\end{equation}
Then the functions
\begin{eqnarray}
\phi_k^\lambda(\xi)=\alpha^{1\over2}e^{i(k+\lambda)\alpha\xi}\phi(\alpha
\xi),~ ~ k\in\mathbb{Z} \label{basis}
\end{eqnarray}
form an orthonormal basis for the Fourier transform
${\cal{F}}V^{\alpha,\lambda}$  of $V^{\alpha,\lambda}$ space. And
hence the functions
$\hat\phi_k^\lambda(x)=\alpha^{-{1\over2}}\hat\phi(\alpha^{-1}(x-k-\lambda)$
form an orthonormal basis of $V^{\alpha,\lambda}$, which give us a family of isometric injections:
\begin{equation}
I_{\alpha,\lambda}: {\cal{V}}^{\alpha,\lambda}\tilde  \rightarrow
V^{\alpha,\lambda}\subset L^2(\mathbb{R})
\end{equation}
The Meyer proposition stipulates that for all $\alpha>0$, and
$\lambda\in\mathbb{R}\slash\mathbb{Z}$ we have:\\

 1. $V^{\alpha,\lambda}\subset V^{\alpha',\lambda'}$ when
 $\alpha(\mathbb{Z}+\lambda) \subset\alpha '(\mathbb{Z}+\lambda')$

 2. $\bigcap_{n\in\mathbb{Z}}V^{2^{n}\alpha, 2^{-n}\lambda} =\left\{0\right\}$, and
  $\overline{\bigcup_{n\in\mathbb{Z}}V^{2^{n}\alpha,2^{-n}\lambda}}=L^2(\mathbb{R})$,

 3. For all $a>0$, $f(x)\in V^{\alpha,\lambda}\Leftrightarrow a^{-1/2}f(a^{-1}x)\in V^{a\alpha,
 \lambda}$\\
 and for all $b\in \alpha\mathbb{Z}$, if  $f(x)\in V^{\alpha,\lambda}$ then $f(x-b)\in V^{\alpha,\lambda+b/\alpha}$
(i.e. the injections $I_{\alpha,\lambda}$ commute  with dilations
and translations). \\

The family of the isometric injections $I_{\alpha,\lambda}$ is
called \textit {multiresolution analysis} of $L^2(\mathbb{R})$ and
is
 $\infty$\textit {-regular} \footnote{A multiresolution analysis is said  {\em $\infty$-regular},
 if any Dirac mass is sent by $I_{\alpha,\lambda}$ on  a function of the Schwartz space $\cal{S}(\mathbb
 R)$.}. It furnishes a precise meaning to the intuitive fact that ${\cal{H^{\alpha,\lambda}}}$ "tends"
 to $L^2(\mathbb{R})$ when $\alpha$ goes to $0$. To construct the
 isometric injections on two copies of $L^2(\mathbb{R})$ it is
 enough  to identify in natural way the space ${\cal{H^{\alpha,\lambda}}}$ to
 two copies of  ${\cal{H}}^{2\alpha,\lambda/2}$. In fact considering
 the operators
\begin{eqnarray}
 {\cal U}_{\alpha,\lambda}:{\cal H}^{\alpha,\lambda}\rightarrow {\cal H}^{\alpha,\lambda}  ,\qquad {\cal J}:{\cal H}^{2\alpha,\lambda/2}\rightarrow {\cal H}^{\alpha,\lambda} \nonumber
\end{eqnarray}
defined by
\begin{equation}
 {\cal U}_{\alpha,\lambda}.f(\psi)=e^{-i\pi\lambda}f\left(\psi+\pi\alpha^{-1}\right), \qquad {\cal J}.f(\psi)=\left(1+e^{i\alpha \psi}\right)f\left(\psi\right).
\end{equation}
The application
\begin{equation}
({\cal J},~ {\cal U}_{\alpha,\lambda}\circ{\cal J}):{\cal
H}^{2\alpha,\lambda/2}\oplus {\cal
H}^{2\alpha,,\lambda/2}\rightarrow {\cal H}^{\alpha,\lambda}
\end{equation}
is an isometric isomorphism with inverse given by
\begin{equation}
\left(
\begin{array}{cc}
  {\cal R} \\
  {\cal R}\circ {\cal U}_{\alpha,\lambda}
\end{array}
\right)
\end{equation}
where ${\cal R}$ is the adjoint of the injection ${\cal J}$ and reads
 \begin{equation}
 {\cal R}.f(\psi)={1\over4}\left(1+e^{-i\alpha \psi}\right)f\left(\psi\right)+
 {1\over4}\left(1-e^{-i\alpha \psi}\right)f\left(\psi+\pi\alpha^{-1}\right).
\end{equation}
Now let $I_{\alpha,\lambda}$  a multiresolution analysis of
Littlewood-Paley-Meyer of  $L^2(\mathbb{R})$, so for all $f \in
{\cal H}^{\alpha,\lambda}$ we define ${\cal I}_{\alpha,\lambda}$ by:
\begin{equation}
{\cal I}_{\alpha,\lambda}(f)(\psi):={\cal F}^{-1}\circ
I_{\alpha,\lambda}\circ {\cal F}(f)(\psi)=\phi(\alpha \psi)f(\psi)
\end{equation}
If we set, for any $u\in L^2(\mathbb{R})$,
\begin{equation}
{\cal A}_{\alpha,\lambda}u(\psi)=\sum_{k\in
\mathbb{Z}}\phi(\alpha\psi+2k\pi)e^{-2ik\pi\lambda}u(\psi+2k\pi
\alpha^{-1})
\end{equation}
then  ${\cal A}_{\alpha,\lambda}u\in{\cal H}^{\alpha,\lambda}$. By
the virtue of (\ref{per}) we have, for any $f\in {\cal
H}^{\alpha,\lambda}$, $({\cal A}_{\alpha,\lambda}\circ {\cal
I}_{\alpha,\lambda})f(\psi)=f(\psi)$.\\
We can check easily that ${\cal A}_{\alpha,\lambda}$ is the adjoint
of  ${\cal I}_{\alpha,\lambda}$ and that  ${\cal
P}_{\alpha,\lambda}= {\cal I}_{\alpha,\lambda}\circ{\cal
A}_{\alpha,\lambda}$ is an orthogonal projector of $L^2(\mathbb{R})$
on the Fourier transform ${\cal F}V^{\alpha,\lambda}$ of the space
 $V^{\alpha,\lambda}$. And then for all $u, ~v\in
\cal{S}(\mathbb{R})$ and  $X \in V$, the following expression
\begin{equation}
\left(
\begin{array}{cc}
 {\cal I}_{2\alpha,\lambda/2}  & 0 \\
   0 & {\cal I}_{2\alpha,\lambda/2}
\end{array}
\right)\left(
\begin{array}{c}
 {\cal R}   \\
    {\cal R}\circ {\cal U}_{\alpha,\lambda}
  \end{array}
\right)\circ R^{\alpha,\lambda}_h(exp_\alpha X) \circ \left({\cal
J},~{\cal U}_{\alpha,\lambda}\circ {\cal J}\right)\left(
\begin{array}{cc}
 {\cal A}_{2\alpha,\lambda/2}  & 0 \\
   0 & {\cal A}_{2\alpha,\lambda/2}
\end{array}
\right) \left(
\begin{array}{c}
 u   \\
    v
  \end{array}
\right)
\end{equation}
tends to
\begin{equation}
\left(
\begin{array}{cc}
 R^{ h}(exp_0 X)  & 0 \\
   0 & R^{ -h}(exp_0 X)
\end{array}
\right) \left(
\begin{array}{c}
 u   \\
    v
  \end{array}
\right)
\end{equation}
when $\alpha$ tends to $0$, the convergence holding in
$\cal{S}(\mathbb{R})\oplus\cal{S}(\mathbb{R})$ in Fr\'echet sense
and uniformally for all $X$ belong to a compact set of $V$.


\newpage
\begin{thebibliography}{99}

\bibitem{AS}
M. Abramowitz and I. Stegun, {\em Handbook of Mathematical
Functions}, Dover, New York, (1964).

\bibitem{AM}
M. Andler et D. Manchon,  Op\'erateurs aux Diff\'erences Finies,
Calcul Pseudo-diff\'erentiel et Repr\'esentations des Groupes de
Lie, Journal of geometry and Physics 27, 1-29 (1998).

\bibitem{Db}
I. Daubechies, {\em Ten lectures on Wavelets}, Philadelphia,
(1992).

\bibitem{Ince}
 E. L. Ince, {\em Ordinary Differential Equations}, Dover,
New York (1957).

\bibitem{KM}
E. G. Kalnins and Willard Miller Jr.,  Lie Theory and Separation
of Variables. 4. The Groupe $SO(2,1)$ and $SO(3)$, J. Math. Phys.
V. 15, No. 8, 1263-1272 (1974).

\bibitem{KMP}
E. G. Kalnins, Willard Miller Jr. and G. S. Pogosyan, Contraction
of Lie Algebras: Application to Special Functions and Separation
of Variables, J. Phys. A: Math. Gen. { \bf 32}, 4709-4732 (1999).

\bibitem{Mal}
S. Mallat,  Multiresolution Approximation and Wavelets Orthonormal
Bases of $L^2(\mathbb{R})$, Trans. Amer. Math. Soc., vol. 315,
69-87, (1989).

\bibitem{Meixner}
J. Meixner and F.W. Schafke, {\em Mathieusche Funktionen und
Spheroidfunktionen}. Springer Berlin, Gottingen, Heidelberg,
(1954).

\bibitem{Me}
Y. Meyer, {\em Ondelettes et Op\'erateurs}, tome 1, Hermann,
Paris, (1990).

\bibitem{Meck}
J. Mickelsson and J. Niederle, Contractions of Representations of
de Sitter Groups, Commun. math. Phys. 27, 167-180 (1972).

\bibitem{MI}
Willard Miller Jr., {\em Symmetry and Separation of Variables},
Addison-Wesley, Reading, massachusettes, (1977).

\bibitem{SLS}
S. Yu. Slavyanov, W. Lay and A. Seeger, {\em Special Fuctions a
Unified Theory Based on Singularities}, Oxford Univ. Press, New York
(2000).

\bibitem{slav1}
S. Yu. Slavyanov, W. Lay, A. M. Akopyan, A. B. Pirozhnikov, V. Yu.
Dmitriev, A. B. Yazik, and V. Zhegunov,  A Knowledge Base On
Special Functions, Journal of Mathematical Sciences, Vol. 108, No.
6, (2002).

\bibitem{TA}
Michael. E. Taylor, { \em Noncommutative Harmonic Analysis}, AMS.
(1986).

\bibitem{Th}
S. Thangavelu, {\em Harmonic Analysis On the Heisenberg Group },
Progress in Mathematics, Birkh\"{a}user V. 159 (1998).

\bibitem{W}
P. Winternitz, I. luk\u{a}c, and Y. Smorodinski\u{i},  Quantum
Numbers in the Little Groups of the Poincar\'e Group, Sov. J.
Nucl. Phys. 7, 139 (1968).

\bibitem{YO}
M. Yamada and K. Ohkitani, An Identification of Energy Cascade in
Turbulence by Orthonormal Wavelet Analysis, Progress of
Theoretical Physics, V. 86, No. 4, (1991).
\end{thebibliography}
\end{document}